\documentclass{amsart}
\usepackage{amsmath, amsthm, amsfonts, amssymb, graphicx, hyperref, verbatim, cleveref, enumitem}

\theoremstyle{plain}

\newtheorem{theorem}{Theorem}[section]
\newtheorem{remark}[theorem]{Remark}

\newtheorem{lemma}[theorem]{Lemma}
\newtheorem{defn}[theorem]{Definition}
\newtheorem{coro}[theorem]{Corollary}
\newtheorem{prop}[theorem]{Proposition}

\newcommand{\R}{\mathbb{R}}

\newcommand{\bB}{\mathbb{B}}

\newcommand{\Bb}{\mathcal{B}_{\operatorname{broad}}}
\newcommand{\Bn}{\mathcal{B}_{\operatorname{narrow}}}

\newcommand{\N}{\mathbb{N}}

\newcommand{\Yb}{Y_{\operatorname{broad}}}
\newcommand{\Yn}{Y_{\operatorname{narrow}}}

\newcommand{\W}{\mathbb{W}}
\renewcommand{\S}{\mathbb{S}}
\newcommand{\T}{\mathbb{T}}
\newcommand{\les}{\lesssim}

\newcommand{\sq}{\square}
\newcommand{\ges}{\gtrsim}
\newcommand{\less}{\lessapprox}

\newcommand{\eps}{\epsilon}

\AtBeginDocument{%
	\def\MR#1{}
}

\begin{document}
	\title{Weighted decoupling estimates and the {B}ochner-{R}iesz means}	
	\author{Jongchon Kim}
	\address{Department of Mathematics, City University of Hong Kong, Hong Kong SAR, China}
	\email{jongchon.kim.work@gmail.com}
	
	\begin{abstract} 
		We prove new weighted decoupling estimates. As an application, we give an improved sufficient condition for almost everywhere convergence of the Bochner-Riesz means of arbitrary $L^p$ functions for $1<p<2$ in dimensions 2 and 3.
	\end{abstract}
	
	\maketitle
\section{Introduction}
Let $d\geq 2$. For $\lambda, t>0$, the Bochner-Riesz means $S^\lambda_t f$ are defined as 
\[  \widehat{S^\lambda_t f}(\xi) =  \Big( 1-\frac{|\xi|^2}{t^2} \Big)_+^\lambda  \widehat{f}(\xi), \; \xi\in \R^d. \]
We consider the pointwise convergence problem of the Bochner-Riesz means; for which $\lambda$, $S^\lambda_t f$ converges to $f$ almost everywhere as $t\to \infty$ for arbitrary $f\in L^p(\R^d)$? For $p> 2$, this is the case if $\lambda>\max(d(\frac{1}{2}-\frac{1}{p})-\frac{1}{2}, 0)$. This result is due to Carbery  \cite{CarM} for $d=2$ and Carbery, Rubio de Francia and Vega  \cite{CRV} for every $d\geq 2$ (see also \cite{LeeSeegerConv} for endpont results).  Note that the given range of $\lambda$ is  precisely the one for which $S^\lambda_t$ is known to be bounded on $L^p(\R^2)$  and is conjectured to be bounded on $L^p(\R^d)$ for $d>2$ (see \cite{CaSj,Fefferman,CorBR} for $n=2$ and \cite{Guo_Bochner} and references therein for partial results for $n\geq 3$). Meanwhile, it is not known what happens when $\lambda =0$ and $p=2$, unlike the one dimensional case where the Carleson-Hunt theorem is available.

For the case  $1\leq p\leq 2$, by Stein's maximal principle \cite{SteinL}, almost everywhere convergence of an arbitrary $L^p$ function is equivalent to the following weak-type estimate
\begin{equation}\label{eqn:maxBR}
	\|\sup_{t>0} | S^\lambda_t f| \|_{L^{p,\infty}(\R^d)} \les \| f\|_{L^p} 
\end{equation}
for test functions $f$. Here, by $A\les B$, we mean $A\leq CB$ for an absolute positive constant $C$.  Examples of Fefferman \cite{Fefferman_ball} and Tao \cite{TaoW} show that  \eqref{eqn:maxBR} fails unless $\lambda > \max(0,d(\frac{1}{p}-\frac{1}{2}) - \frac{1}{2p})$. We prove a new partial result on this problem in dimension $d=2$. 
\begin{theorem}\label{thm:BR}
	Let $d=2$ and $p=\frac{86}{57}$. Then  \eqref{eqn:maxBR} holds for any $\lambda>\frac{9}{86}$. 
\end{theorem}

By an interpolation of classical $L^1$ and $L^2$ estimates,  \eqref{eqn:maxBR} holds if $\lambda > (d-1)(\frac{1}{p}-\frac{1}{2})$. For $d=2$, Tao \cite{TaoM} improved this classical sufficient condition by giving a $L^{10/7}$-estimate. Building on \cite{TaoM}, Li and Wu  \cite{LiWu}  gave a new $L^{18/13}$-estimate using the Bourgain-Demeter decoupling theorem \cite{BD}.  More recently, Gan and Wu \cite{GanWu} obtained an improved $L^{10/7}$-estimate by using a weighted version of the $\ell^p$-decoupling theorem.  Moreover, they gave a new partial result in dimension $d=3$. See \Cref{thm:BR3d} for an improved bound in $\R^3$.

Combining \Cref{thm:BR} with the classical $L^1$ and $L^2$ estimates and the $L^{18/13}$-estimate \cite{LiWu} yields  the following corollary, which improves previously known results for  $18/13<p<2$.
\begin{coro}\label{cor:BR}
	Let $d=2$, $1< p< 2$ and \[ \lambda>\max\left(\frac{24}{23p}-\frac{27}{46}, \frac{9}{14p}-\frac{9}{28}, \frac{6}{5p}-\frac{7}{10} \right).\] 
	Then \eqref{eqn:maxBR} holds. Consequently, $\lim_{t\to \infty} S^\lambda_t f(x) = f(x)$ for almost every $x\in \R^2$ for any $f\in L^p(\R^2)$.
\end{coro}

Our new estimate is based on weighted $\ell^p$-decoupling inequalities to be discussed below. We first recall the $\ell^p$-decoupling theorem which follows from the Bourgain-Demeter $\ell^2$-decoupling theorem \cite{BD}. Let $S\subset \R^d$ be a strictly convex $C^2$ hypersurface with Gaussian curvature comparable to $1$. For $R\gg 1$, cover the $R^{-1}$-neighborhood of $S$ by finitely overlapping rectangles $\theta$ of dimensions $R^{-1/2} \times \cdots \times R^{-1/2} \times R^{-1}$. Suppose that the Fourier transform of $f_\theta$ is supported on $\theta$ and $f=\sum_\theta f_\theta$. The $\ell^p$-decoupling inequality says that, for $2\leq p\leq \frac{2(d+1)}{d-1}$ and any ball $B_R\subset \R^d$ of radius $R$,  
\[ 	\| f \|_{L^p(B_R)} \less  R^{\frac{d-1}{2} (\frac{1}{2}-\frac{1}{p})}  \Big( \sum_{\theta}  \| f_\theta \|_{L^p}^p \Big)^{\frac{1}{p}}, \]
where we denote by $A\less B$ the inequality $A \leq C_\eps R^\eps B$ which holds for any $\eps>0$.

Gan and Wu \cite{GanWu} observed that the $\ell^p$-decoupling estimate can be improved when the integration of $f$ is taken over a subset $Y\subset B_R$ for the exponent $p=\frac{2d}{d-1}$ and showed that it yields an improved weak-type estimate \eqref{eqn:maxBR}. Here we give a refinement of their weighted $\ell^p$-decoupling estimates for a wider range of $p$.
	\begin{theorem}\label{thm:main0} Let  $p_d=\frac{2(d+1)}{d-1}$,  $2\leq p \leq p_d$, and
		\begin{equation}\label{eqn:alpha}
			\alpha(p) = \min\left(\frac{d+1}{2d}\Big(\frac{1}{p}-\frac{1}{p_d}\Big),  \frac{1}{d-1}\Big(\frac{1}{2}-\frac{1}{p}\Big) \right).
		\end{equation}
		Then for any $Y\subset B_R$,
	\begin{align}\label{eqn:main0}
		\| f \|_{L^p(Y)} \less (|Y| R^{-d})^{\alpha(p)} R^{\frac{d-1}{2} (\frac{1}{2}-\frac{1}{p})}   \Big( \sum_{\theta} \| f_\theta \|_{L^p}^p \Big)^{\frac{1}{p}}.
			\end{align}
\end{theorem}

	Note that $\alpha(p)>0$ for $2<p<p_d$, so the  factor $(|Y|R^{-d})^{\alpha(p)}$ represents a gain over the $\ell^p$-decoupling inequality. Gan and Wu \cite{GanWu} proved \eqref{eqn:main0} for the exponent $p=\frac{2d}{d-1}$ with $\alpha(p)=\frac{d+1}{2d}\Big(\frac{1}{p}-\frac{1}{p_d}\Big) = \frac{1}{16}$ for $d=2$ and with $\alpha(p) =  \frac{1}{3d-1}\Big(\frac{1}{2}-\frac{1}{p}\Big) = \frac{1}{2d(3d-1)}$ for $d\geq 3$. 	  
	
	We next discuss sharpness of \Cref{thm:main0}. The exponent $\alpha(p)$ given in \eqref{eqn:alpha} is optimal when  $\frac{2(d^2+2d-1)}{d^2+1} \leq p\leq p_d$, which is the range where $\alpha(p)=\frac{d+1}{2d}(\frac{1}{p}-\frac{1}{p_d})$. Indeed, 
 if \eqref{eqn:main0} holds for any $Y\subset B_R$, then $\alpha(p)$ must obey
 \begin{equation}\label{eqn:alphabound}
 	\alpha(p) \leq \min\left(\frac{d+1}{2d}\Big(\frac{1}{p}-\frac{1}{p_d}\Big),  \frac{1}{2}-\frac{1}{p} \right).
 \end{equation}
	To see this, one takes $\widehat{f_\theta} = \frac{1}{|\theta|}\chi_{\theta}$ for a bump function $\chi_\theta$ supported on $\theta$ and let $Y=B_c(0)$ for sufficiently small $c\sim 1$ (cf. \cite{GanWu}). Then $|\sum_\theta f_\theta(x)| \sim  \# \theta \sim R^{\frac{d-1}{2}}$ on $Y$. This example gives $\alpha(p)\leq \frac{d+1}{2d}\left(\frac{1}{p}-\frac{1}{p_d}\right)$. To see the other upper bound, it suffices to take $f=f_\theta$ for a fixed $\theta$ and take $Y=\theta^*$, the rectangle dual to $\theta$. 
	
	Comparing \eqref{eqn:alphabound} and \eqref{eqn:alpha}, we see that the exponent $\alpha(p)$ given in \Cref{thm:main0} is optimal for $d=2$ for all $2\leq p\leq 6$. Note that $\alpha(p)$ takes the maximum value $\alpha(p)=\frac{1}{7}$ at $p=\frac{14}{5}$; we will use \Cref{thm:main0} with this particular $p$ for our application to the Bochner-Riesz means. On the other hand, in dimensions  $d\geq 3$, we do not currently have any reason to believe that the exponent 	$\alpha(p) =  \frac{1}{d-1}\big(\frac{1}{2}-\frac{1}{p}\big)$ from \eqref{eqn:alpha} for $p$ close to $2$  is optimal. Moreover, if we assume further conditions on $Y\subset B_R$, then it is possible that \eqref{eqn:main0} may hold with  $\alpha(p)$ which does not obey  \eqref{eqn:alphabound}; see \cite{GanWu} for a discussion of such a possibility when $|Y|\sim R^{1/2}$ and $d=2$.
	
	Our proof of \Cref{thm:main0} builds on \cite{GanWu}, which is based on a broad-narrow analysis going back to the work of Bourgain and Guth \cite{BG} (see also \cite{Guth2}). In the broad case, we use the Bennett-Carbery-Tao  multilinear restriction estimate \cite{BCT} to exploit the transversality and combine it with the $\ell^p$-decoupling theorem as in \cite{GanWu}. In the narrow case, we do induction on scales in a way similar to the one developed in the proof of the refined Strichartz estimate by Du, Guth and Li \cite{DGL}. The novelty of our proof lies on the analysis of the narrow case; we utilize refined and weighted decoupling estimates due to Guth, Iosevich, Ou and Wang \cite{GIOW} and Du, Ou, Ren and Zhang  \cite{du2023weighted} in order to perform induction on scales more efficiently. Given  \Cref{thm:main0}, \Cref{thm:BR} follows from the argument of \cite{GanWu}. 
		
	The paper is organized as follows. In \Cref{sec:weighted}, we prove a multilinear weighted decoupling estimate which will be used in the analysis of the broad part. We also state a weighted refined decoupling estimate to be used in the narrow part. In \Cref{sec:main}, we prove \Cref{thm:main0}. In \Cref{sec:BR}, we give an exposition of the paper \cite{GanWu} on the implication of weighted $\ell^p$-decoupling estimates to the weak-type estimate \eqref{eqn:maxBR} for the Bochner-Riesz means, completing the proof of \Cref{thm:BR}. In addition, we give an improved bound in dimension $d=3$. Throughout the paper, we denote $p_d = \frac{2(d+1)}{d-1}$ with the convention $p_1 = \infty$.
	
\section{Weighted and refined decoupling inequalities} \label{sec:weighted}
In this section, we prove a multilinear weighted refined decoupling inequality, \Cref{thm:refinedGW}, and then introduce a weighted refined decoupling estimate, \Cref{thm:b}, which will be used in the proof of \Cref{thm:main0}. 

\subsection{A multilinear weighted decoupling estimate}
We first recall the multilinear restriction estimate \cite{BCT}. Suppose that $U_j\subset \S^{d-1}$ for $j=1,2,\cdots,d$, such that for any $v_j\in U_j$, 
\begin{equation}\label{eqn:transversality}
	|v_1 \wedge v_2 \wedge \cdots \wedge v_d| \geq \nu,
\end{equation}
for some $0<\nu\leq  1$. 
\begin{theorem}[Multilinear restriction]\label{thm:multiRest} 
	For each $1\leq j\leq d$, let $S_j \subset S$ be a cap such that the set $U_j\subset \S^{n-1}$ of vectors normal to $S_j$ satisfies \eqref{eqn:transversality}.
	Suppose that the Fourier transform of $f_j$ is supported on the $R^{-1}$-neighborhood of $S_j$. Then 
	\[ \| \prod_{j=1}^d |f_j|^{1/d}  \|_{L^{\frac{2d}{d-1}}(B_R)} \les \nu^{-O(1)}  R^{\epsilon} R^{-1/2}   \prod_{j=1}^d \|  f_j \|_{L^2}^{1/d}.\]
\end{theorem}

Next we recall the refined decoupling estimate which is stated in terms of wave packets. We recall the setup.  Given $0<\eps \ll 1$, let $\eps_1 = \eps^{100}$. For each $\theta$, let $\T_\theta$ denote a collection of finitely overlapping tubes $T$ covering $\R^d$ of length $\sim R^{1+\eps_1}$ and radius $\sim R^{(1+\eps_1)/2}$ with long axis orthogonal to $\theta$. Let $\T(R) =  \cup_\theta \T_\theta$. Given $T\in \T(R)$, let $\theta(T)$ denote the $\theta$ for which $T\in \T_\theta$. 

We say that a function $f_T$ is microlocalized to $(T,\theta)$ if $f_T$ and $\widehat{f_T}$ are essentially supported on $2T$ and $2\theta$, respectively, in the sense that the $L^p$ norms of the restrictions of $f_T$ to $(2T)^c$ in the physical space and to $(2\theta)^c$ in the frequency space are both $O(R^{-100d} \| f_T \|_{L^p})$. Terms involving $O(R^{-100d})$ can be absorbed into the main term of the estimate, so we will ignore these terms in the following statements and proofs.  We recall that a function $f$ whose Fourier transform is supported on the $R^{-1}$-neighborhood of $S$ admits a wave packet decomposition $f=\sum_\theta f_\theta =\sum_{T\in \T(R)} f_T$, where $f_T$ is  microlocalized to $(T,\theta(T))$.

We can now state the refined decoupling theorem, \cite[Theorem 4.2]{GIOW}.
\begin{theorem}[Refined decoupling] \label{thm:refined_decoupling}
Let $\W \subset \T(R)$ and suppose that $f=\sum_{T\in \W} f_T$, where $f_T$ is microlocalized to $(T,\theta(T))$. 
	Let $Y\subset B_R$ be the union of a collection of $R^{1/2}$-cubes $Q$, each of which is contained in $3T$ for at most $M\geq 1$ tubes $T\in \W$. Then for $2\leq p\leq p_d$,
	\[ \| f \|_{L^p(Y)} \less M^{\frac{1}{2}-\frac{1}{p}} \left( \sum_{T\in \W} \| f_T \|_{L^p}^p \right)^{1/p}. \]
\end{theorem}

The following is a  multilinear refinement of \Cref{thm:refined_decoupling}.
\begin{theorem}\label{thm:refinedGW}
	For each $j=1,2,\cdots, d$, let $S_j \subset S$ be a cap such that the set $U_j\subset \S^{n-1}$ of vectors normal to $S_j$ satisfies \eqref{eqn:transversality}. Let $\W_j \subset \T(R)$ such that $2\theta(T)$ is contained in the $O(R^{-1})$-neighborhood  of $S_j$ for every $T\in \W_j$ and let $f_j = \sum_{T\in \W_j} f_T$, where $f_T$ is microlocalized to $(T,\theta(T))$. Let $\{Q\}$ be a collection of $R^{1/2}$-cubes contained in $B_R$ such that each cube in the collection is contained in $3T$ for at most $M_j$ tubes $T\in \W_j$ for each $1\leq j\leq d$. Let $Y$ be a subset of  $\cup_Q Q$. Then for $2 \leq p\leq p_d$, 
	\begin{align*}
	\left\| \prod_{j=1}^d |f_j|^{\frac{1}{d}} \right\|_{L^{p}(Y)}  \less  \nu^{-O(1)} ( |Y| R^{-d})^{\frac{d+1}{2d} (\frac{1}{p} - \frac{1}{p_d})} \prod_{j=1}^d \left[ M_j^{\frac{1}{2}-\frac{1}{p}} \Big( \sum_{T\in \W_j} \| f_T \|_{L^p}^p \Big)^{\frac{1}{p}} \right]^{\frac{1}{d}}.
\end{align*}
\end{theorem}
For the proof, we interpolate \Cref{thm:multiRest}  and \Cref{thm:refined_decoupling}. We remark that the use of \emph{refined} decoupling is optional here as we shall use \Cref{thm:refinedGW} with $M_j \les R^{\frac{d-1}{2}}$.  We note that \Cref{thm:refinedGW} for the case $p=\frac{2d}{d-1}$ was essentially proved in \cite{GanWu}. Our proof is slightly more direct than the one given there in that it doesn't go through the multilinear Kakeya estimate. 

\begin{proof}
	We prove the result for $2\leq p\leq \frac{2d}{d-1}$; the result for the remaining $p$ follows from interpolation with the statement for $p=p_d$, which is a consequence of H\"{o}lder and \Cref{thm:refined_decoupling}. 	
	
We may view $f_T$ as essentially constant on $T$; this can be made precise by using the frequency localization property of $f_T$, cf. \Cref{sec:broad}. By dyadic pigeonholing, we may reduce the matter to the case where, for each $j=1,\cdots,d$,  $| f_T | \sim \gamma_j$  on $T$ for some $\gamma_j>0$ for all $T\in \W_j$. By homogeneity, we may assume that $\gamma_j=1$. Then for $1\leq q< \infty$, 
\[ \sum_{T\in \W_j} \| f_T \|_{L^q}^q  \sim  |T|  |\W_j |.\] 
	
By H\"{o}lder, we get
	\begin{equation}\label{eqn:roughlyconst}
		\begin{split}
			\left\| \prod_{j=1}^d |f_j|^{\frac{1}{d}} \right\|_{L^{p}(Y)} \les |Y|^{\frac{1}{p}-\frac{d-1}{2d}} 	\left\| \prod_{j=1}^d |f_j|^{\frac{1}{d}} \right\|_{L^{\frac{2d}{d-1}}(Y)}. 
		\end{split}
	\end{equation}  	
The multilinear restriction estimate and the $L^2$-orthogonality imply
	\begin{equation}\label{eqn:muld}
 	\left\| \prod_{j=1}^d |f_j|^{\frac{1}{d}} \right\|_{L^{\frac{2d}{d-1}}(B_R)} \less R^{-\frac{1}{2}} \prod_{j=1}^d \| f_j \|_{L^2}^{\frac{1}{d}} \less R^{-\frac{1}{2}} \prod_{j=1}^d (|T||\W_j| )^{\frac{1}{2d}}. 	
	\end{equation}
	
	 On the other hand, by H\"{o}lder and the refined decoupling estimate, \Cref{thm:refined_decoupling}, we have
\begin{equation}\label{eqn:red}
	\begin{split}
		\left\| \prod_{j=1}^d |f_j|^{\frac{1}{d}} \right\|_{L^{\frac{2d}{d-1}}(Y)} &\les  |Y|^{\frac{d-1}{2d}-\frac{1}{p_d}} \left\| \prod_{j=1}^d |f_j|^{\frac{1}{d}} \right\|_{L^{p_d}(Y)}  \\
	& \less |Y|^{\frac{d-1}{2d}-\frac{1}{p_d}} \prod_{j=1}^d \left[ M_j^{\frac{1}{2}-\frac{1}{p_d}} \left( \sum_{T\in \W_j} \| f_T \|_{L^{p_d}}^{p_d}  \right)^{\frac{1}{p_d}}  \right]^{\frac{1}{d}} \\
	 & \less |Y|^{\frac{d-1}{2d(d+1)}} \prod_{j=1}^d \left[ M_j^{\frac{1}{d+1}} (|T| |\W_j|)^{\frac{1}{p_d}} \right]^{\frac{1}{d}}.
		\end{split}
\end{equation}
We combine the $(d+1)(\frac{1}{p}-\frac{1}{p_d})$-th power of \eqref{eqn:muld} and $(d+1)(\frac{1}{2}-\frac{1}{p})$-th power of \eqref{eqn:red}. Then we get 
	\[ 	\left\| \prod_{j=1}^d |f_j|^{\frac{1}{d}} \right\|_{L^{\frac{2d}{d-1}}(Y)} \less R^{-\frac{d+1}{2} (\frac{1}{p}-\frac{1}{p_d})}  |Y|^{\frac{d-1}{2d}(\frac{1}{2}-\frac{1}{p})} \prod_{j=1}^d \left[ M_j^{\frac{1}{2}-\frac{1}{p}} \left( \sum_{T\in \W_j} \| f_T \|_{L^p}^{p}  \right)^{\frac{1}{p}}  \right]^{\frac{1}{d}}.  \]
Finally, we combine this estimate with \eqref{eqn:roughlyconst}. 
\end{proof}

\subsection{A weighted refined decoupling estimate}\label{sec:weird}
Let $v(T) \in \S^{d-1}$ denote the direction of $T\in \T$. Following \cite{du2023weighted}, for $R^{-1/2}\leq r\leq 1$ and $1\leq m\leq d$, we say that $f = \sum_{T\in \W} f_T$ has $(r,m)$-concentrated frequencies if there exists a  $m$-dimensional subspace $V\subset \R^d$ such that $\angle(v(T), V) \leq r$ for any $T \in \W$. 

When $f$ has $(R^{-1/2},m)$-concentrated frequencies, one can decouple $f$ by $\ell^2$-decoupling in dimension $m$ (cf. \cite[Lemma 9.3]{Guth2}). In the paper  \cite{du2023weighted}, it was observed that this continues to hold for the refined decoupling estimate as well.
\begin{theorem}[{\cite[Theorem 1.2(a)]{du2023weighted}}\footnote{The case $m=1$ is omitted in \cite{du2023weighted}, but it trivially holds as  $O(1)$ boxes $\theta$ are involved.}] \label{thm:refined_decouplingDu}
	Let $\W \subset \T(R)$ and suppose that $f=\sum_{T\in \W} f_T$, where $f_T$ is microlocalized to $(T,\theta(T))$. 	Let $Y\subset B_R$ be the union of a collection of $R^{1/2}$-cubes $Q$, each of which is contained in $3T$ for at most $M\geq 1$ tubes $T\in \W$. Suppose that $f$ has $(R^{-1/2},m)$ concentrated frequencies for some $1\leq m\leq d$.  Then for $2\leq p\leq p_m$,
	\[ \| f \|_{L^p(Y)} \less M^{\frac{1}{2}-\frac{1}{p}} \left( \sum_{T\in \W} \| f_T \|_{L^p}^p \right)^{1/p}. \]
\end{theorem}

A consequence of \Cref{thm:refined_decouplingDu} is a  weighted refined decoupling estimate with weights $w: Y\to [0,1]$ such that $\int_{Q} w \les R^{\alpha/2}$ holds for any $R^{1/2}$-cube $Q$ in $Y$; see \cite[Theorem 1.2(b)]{du2023weighted}). Here we state a slightly more general version which does not require such an assumption. 
\begin{theorem}[Weighted refined decoupling] \label{thm:b}
	Let $1\leq m\leq d$ and $2\leq p\leq p_m$. Suppose that $f= \sum_{T\in \W} f_T$, where $\W \subset \T(R)$ and each $f_T$ is microlocalized to $(T,\theta(T))$. Assume that $f$ has $(R^{-1/2},m)$-concentrated frequencies. Let $\{ Q\}$ be a collection of disjoint $R^{1/2}$-cubes in $B_R$ such that each $Q$ is contained in $3T$ for at most $M\geq 1$ tubes $T\in \W$. Let $Y$ be a subset of  $\cup_Q Q$. Let $w: Y \to [0,1]$ and $w(T) := \int_T w$. Then
	\[  \| f\|_{L^p(w)}  \less \left( \frac{w(3T)}{|T|} \right)^{\frac{1}{p}-\frac{1}{p_m}}  M^{\frac{1}{2}-\frac{1}{p}}   \left( \sum_{T\in \W} \| f_T\|_{L^p}^p \right)^{1/p}.\]
\end{theorem}
The power of $\frac{w(3T)}{|T|}$  represents a gain over the refined decoupling estimate. \Cref{thm:b} essentially follows from the proof of \cite[Theorem 1.2(b)]{du2023weighted}). We include a sketch of the proof for the sake of completeness.
\begin{proof} 
	By dyadic pigeonholing, we may assume that each $R^{1/2}$-cube $Q$ in the collection is contained in $3T$ for $\sim M$  tubes $T\in \W$. By H\"{o}lder,
\[ \| f\|_{L^p(w)} \les w(Y)^{\frac{1}{p}-\frac{1}{p_m}} \| f\|_{L^{p_m}(Y)}.\]	
	Next, we use  \Cref{thm:refined_decouplingDu} with  exponent $p=p_m$ and then replace the $\ell^{p_m}(L^{p_m})$ norm back to $\ell^{p}(L^{p})$ as in the proof of \Cref{thm:refinedGW}, which give
	\begin{eqnarray}\label{eqn:weidec}
		\| f\|_{L^p(w)} \less \left(\frac{w(Y) M}{|\W| |T| }\right)^{\frac{1}{p}-\frac{1}{p_m}} M^{\frac{1}{2}-\frac{1}{p}}   \left( \sum_{T\in \W} \| f_T\|_{L^p}^p \right)^{1/p}.
	\end{eqnarray}
 Then 
	\[	w(Y) M \sim   \sum_{Q} \sum_{T\in \W: Q \subset 3T} w(Q) =  \sum_{T\in \W} \sum_{Q:Q \subset 3T} w(Q) \leq |\W| \max_{T\in \W} w(3T). \]
Combining this inequality and \eqref{eqn:weidec} finishes the proof.
\end{proof}

\section{Broad-narrow analysis: Proof of \Cref{thm:main0}}\label{sec:main}
In this section, we prove the following theorem which implies  \Cref{thm:main0} by the wave packet decomposition.
\begin{theorem}\label{thm:main} Suppose that $\W \subset \T(R)$ and $f=\sum_{T\in \W} f_T$, where $f_T$ is microlocalized to $(T,\theta(T))$. Let $p_d = \frac{2(d+1)}{d-1}$ and $Y\subset B_R$. Then for $2\leq p \leq p_d$ and $\alpha(p)$ defined in \eqref{eqn:alpha}, 
	\begin{align}\label{eqn:goal}
		\| f \|_{L^p(Y)} \less (|Y| R^{-d})^{\alpha(p)} R^{\frac{d-1}{2} (\frac{1}{2}-\frac{1}{p})}   \left( \sum_{T\in \W} \| f_T \|_{L^p}^p \right)^{\frac{1}{p}}.
	\end{align}
\end{theorem}

We turn to the proof of \Cref{thm:main}. We decompose $f$ at a scale dictated by the parameter $K=R^{\eps_1}=R^{\eps^{100}}$. Consider a cover of $\mathcal{S}$ by rectangles $\tau$ of dimensions $K^{-1}\times \cdots \times K^{-1}\times K^{-2}$ such that for each $T\in \W$, there exists at least one and at most $O(1)$ rectangles $\tau$ containing $\theta(T)$. We fix one such $\tau=\tau(T)$ for each $T\in \W$. We let $\W_\tau= \{ T\in \W: \tau(T) = \tau\}$ and $f_\tau = \sum_{\T \in \W_\tau} f_T$, so that $f=\sum_\tau f_\tau$. 

For each lattice $K^2$-cube $B$ intersecting $Y$, define the ``significant set" 
\[ S(B) = \{ \tau :  \| f_\tau \|_{L^p(B\cap Y)} \geq \frac{1}{100 \#\tau} \| f\|_{L^p(B\cap Y)} \}.  \]
Then the definition of $S(B)$ ensures that for each $K^2$-cube $B$,
\[ \| f\|_{L^p(B\cap Y)} \sim \| \sum_{\tau\in S(B)} f_\tau \|_{L^p(B\cap Y)}. \]

We say that a lattice $K^2$-cube $B$  intersecting $Y$ is narrow if there exists a $(d-1)$-dimensional subspace $V\subset \R^d$ such that for every $\tau\in S(B)$, 
\[ \angle(G(\tau),V) \leq K^{-1},\]
where $G(\tau)\in \S^{d-1}$ denote the direction normal to  $\tau$. Otherwise, we say that $B$ is broad. If $B$ is broad, then there exist $\tau_1,\tau_2,\cdots,\tau_d \in S(B)$ such that
\begin{equation}\label{eqn:transverse}
	|G(\tau_1) \wedge G(\tau_2) \wedge \cdots \wedge G(\tau_d)| \ges K^{-(d-1)}.
\end{equation}

Let $\Bb$ denote the collection of broad $K^2$-cubes and $\Bn$ denote the collection of narrow $K^2$-cubes. We let 
\[ \Yb = \cup_{B\in \Bb} B\cap Y, \;\; \Yn = \cup_{B\in \Bn} B\cap Y. \]
We say that we are in the broad case if $\| f\|_{L^p(\Yb)} \geq \| f\|_{L^p(\Yn)}$. Otherwise, we  say that we are in the narrow case. We handle each case in the following subsections.

\subsection{Broad case}\label{sec:broad}
In this subsection, we denote by $A \less B$ expressions of the form $A\leq K^{O(1)} B$; the loss of $K^{O(1)}$ is harmless as long as it is at most  $R^{\eps/2}$.

In the broad case,  we have $\| f\|_{L^p(Y)} \les \| f\|_{L^p(\Yb)}$.  
Let $B\in \Bb$. Then there exist $\tau_1,\tau_2,\cdots,\tau_d \in S(B)$ satisfying \eqref{eqn:transverse}. We fix such a $d$-tuple and denote $\bar{\tau}(B) = (\tau_1,\cdots,\tau_d)$. Let 
\[ \Gamma = \{ \bar{\tau}(B) : B \in \Bb \}. \]
Since $\# \Gamma \les K^{O(1)}$, by dyadic pigeonholing, there exist $\bar{\tau}=(\tau_1,\cdots,\tau_d) \in \Gamma$ and a sub-collection of broad $K^2$-cubes $\Bb^1$ such that $\bar{\tau}(B) = \bar{\tau}$ for all $B\in \Bb^1$ and 
\[   \| f\|_{L^p(\Yb)}\less  \| f\|_{L^p(\Yb^1)}, \text{ where } \; \Yb^1 = \cup_{B\in \Bb^1} B\cap Y. \]

Let $B\in \Bb^1$. Since $\tau_j \in S(B)$ for each $j$, we have 
\begin{align*}
 	\| f \|_{L^p(B\cap Y)}^p \less  \prod_{j=1}^d \| f_{\tau_j} \|_{L^p(B\cap Y)}^{\frac{p}{d}}.
\end{align*} 
Using the Fourier support property of $f_j$, we morally have the following reverse H\"{o}lder inequality (cf. \cite{DZ}):
\begin{align}\label{eqn:reverseH}
  \prod_{j=1}^d \| f_{\tau_j} \|_{L^p(B\cap Y)}^{\frac{p}{d}} \less   \| \prod_{j=1}^d |f_{\tau_j}|^{1/d} \|_{L^p(B\cap Y)}^{p}.
\end{align} 
We shall first assume  \eqref{eqn:reverseH}, and then make it precise later. Summing \eqref{eqn:reverseH} over $B\in \Bb^1$, we get
\[   \| f\|_{L^p(\Yb)}\less  \| \prod_{j=1}^d |f_{\tau_j}|^{1/d} \|_{L^p(\Yb^1)}. \]

By using \Cref{thm:refinedGW} with the bound $M_j\les R^{\frac{d-1}{2}}$ and $\W_{\tau_j} \subset \W$, we obtain
\begin{equation} \label{eqn:refinedGW}
\begin{split}	 \| f\|_{L^p(\Yb)} &\less (|Y| R^{-d})^{\frac{d+1}{2d}(\frac{1}{p}-\frac{1}{p_d})} \prod_{j=1}^d \left[ M_j^{\frac{1}{2}-\frac{1}{p}} \left( \sum_{T\in \W_{\tau_j}} \| f_T \|_{L^p}^p \right)^{1/p} \right]^{1/d} \\ 
	&\less (|Y| R^{-d})^{\frac{d+1}{2d}(\frac{1}{p}-\frac{1}{p_d})} R^{\frac{d-1}{2} (\frac{1}{2}-\frac{1}{p})}    \left( \sum_{T\in \W} \| f_T \|_{L^p}^p \right)^{1/p},\end{split}
\end{equation}
which implies the claimed estimate \eqref{eqn:goal} as $\alpha(p) \leq \frac{d+1}{2d}(\frac{1}{p}-\frac{1}{p_d})$.

Now we make  the step \eqref{eqn:reverseH}  rigorous following  \cite{hickman2023pointwise}. Since $f_{\tau_j}$ has compact Fourier support, we may write $f_{\tau_j} = f_{\tau_j}* \psi$ for some $\psi$ with compact Fourier support. Therefore, we have $|f_{\tau_j}| \les |f_{\tau_j}|* \Psi_{1}$, where $\Psi_{\rho}(x) = \rho (\rho+|x|)^{-(d+1)}$. Moreover, 
\begin{equation}\label{eqn:locallyconstant}
	 |f_{\tau_j}|^{r} \les_{r}  |f_{\tau_j}|^r *\Psi_{1}
\end{equation}
for all $r>0$. This follows from H\"{o}lder for $r>1$. The case $0<r<1$ can be proved by using Bernstein's inequality; see  \cite[Lemma 5.9]{hickman2023pointwise}. We also note that  $\Psi_{1}(x) = K^{2d}\Psi_{K^{2}}(K^2 x) \less \Psi_{K^2}(x)$. Using \eqref{eqn:locallyconstant} with $r=p/d$, we get
\begin{equation}\label{eqn:locallyconstant2}
	 |f_{\tau_j}|^{\frac{p}{d}} \less |f_{\tau_j}|^{\frac{p}{d}}* \Psi_{K^2}.
\end{equation}

We are going to use the following two properties of $\Psi_{K^2}$; its $L^1$ norm is comparable to $1$  and it is locally constant on balls of radius $K^2$. The latter implies that $|f_{\tau_j}|^{\frac{p}{d}}* \Psi_{K^2}(x) \sim  |f_{\tau_j}|^{\frac{p}{d}}* \Psi_{K^2}(y)$ whenever $|x-y|\les K^2$. Therefore, by \eqref{eqn:locallyconstant2},
\begin{align*}
\prod_{j=1}^d \left( \int_{B\cap Y} |f_{\tau_j}|^{p}(x) dx\right)^{1/d} &\less 
	 \prod_{j=1}^d \left( \int_{B\cap Y} \big( |f_{\tau_j}|^{\frac{p}{d}} * \Psi_{K^2}(x) \big)^d dx\right)^{1/d} \\
&\sim |B\cap Y| \prod_{j=1}^d |f_{\tau_j}|^{\frac{p}{d}} * \Psi_{K^2}(y), \text{ for every } y\in B \\
&\sim \int_{B\cap Y} \prod_{j=1}^d |f_{\tau_j}|^{\frac{p}{d}} * \Psi_{K^2}(x) dx.
\end{align*}
Let $f_{\tau_j,y_j} (x)= f_{\tau_j}(x-y_j)$. Summing this estimate over $B\in \Bb^1$, we get
\begin{align*}
	  \| f\|_{L^p(\Yb)}^p &\less  \int_{Y} \prod_{j=1}^d |f_{\tau_j}|^{\frac{p}{d}} * \Psi_{K^2}(x) dx\\
	  &=  \int_{(\R^d)^d}  \| \prod_{j=1}^d |f_{\tau_j,y_j}|^{1/d} \|_{L^p(Y)}^{p}  \prod_{j=1}^d \Psi_{K^2}(y_j) dy_1 \cdots dy_d.
\end{align*}
By \Cref{thm:refinedGW}, we obtain \eqref{eqn:refinedGW}.

\subsection{Narrow case}
In the narrow case, we  use induction on the scale $R$. In this subsection, we denote by $A \less B$ expressions of the form $A\leq (\log R )^{O(1)} B$ and keep track of powers of $K$. 

We start by considering the base case $R\sim 1$. Since $|\W|\sim 1$, we have 
\[ \| f\|_{L^p(Y)}\les  |Y|^{1/p} \max_{T\in \W}  \|f_T\|_{L^\infty} \les |Y|^{1/p} \max_{T\in \W}   \|f_T\|_{L^p}, \] 
which is better than the claimed estimate since $\alpha(p) \leq 1/p$ and $|Y|\les 1$.

Let $B\in \Bn$. Recall that $\| f\|_{L^p(B\cap Y)} \sim \| \sum_{\tau\in S(B)} f_\tau \|_{L^p(B\cap Y)}$. For each $\tau$, we use a wave packet decomposition:
\[ f_\tau = \sum_{T_1 \in \T_{\tau}} f_{T_1}, \]
where $f_{T_1}$ is microlocalized to $(T_1,\tau)$. Note that $\T_{\tau} \subset \T(K^2)$ and each $T_1\in \T_{\tau}$ has length $\sim K^{2(1+\eps_1)}$ and radius $\sim K^{(1+\eps_1)}$. For dyadic $R^{-1000d} \leq \eta \leq |T_1|$, let $\T_{\tau,\eta}$ denote the collection of $T_1\in \T_{\tau}$ such that $|3T_1 \cap Y| \sim \eta$. The contribution of tubes $T_1 \in \T_\tau$ such that $|3T_1 \cap Y| \les R^{-1000d}$ is negligible; it can be absorbed into the main term by crude estimates. By dyadic pigeonholing, there exists dyadic $\eta$ such that 
\[ \sum_{B\in \Bn} \| f\|_{L^p(B\cap Y)}^p \less \sum_{B\in \Bn}\| \sum_{\tau\in S(B)} \sum_{T_1 \in \T_{\tau,\eta}} f_{T_1} \|_{L^p(B\cap Y)}^p. \]
We fix this $\eta$ and let $\T(K^2;B)$ denote the collection of $T_1\in \cup_{\tau\in S(B)} \T_{\tau,\eta}$ such that $2T_1$ intersects $B$. Then 
\[\| f\|_{L^p(Y)}^p \less \sum_{B\in \Bn}\| \sum_{T_1 \in \T(K^2;B)} f_{T_1} \|_{L^p(B\cap Y)}^p. \]
 
For each  $B\in \Bn$, $\sum_{T_1 \in \T(K^2;B)} f_{T_1}$ has $(K^{-1},d-1)$-concentrated frequencies and $|S(B)| \les K^{d-2}$. Therefore by the weighted refined decoupling  \Cref{thm:b}, 
\begin{equation}\label{eqn:mainb}
	\begin{split}
		\|  \sum_{T_1 \in \T(K^2;B)} f_{T_1}\|_{L^p(B\cap Y)} \les &K^{\eps_1}  (\eta K^{-(d+1)})^{\frac{1}{p}-\frac{1}{p_{d-1}}} \\ &K^{(d-2)(\frac{1}{2}-\frac{1}{p})}
		\left(   \sum_{T_1 \in \T(K^2;B)} \|  f_{T_1} \|_{L^p}^p \right)^{\frac{1}{p}}.
	\end{split}
\end{equation}

For each $\tau$, we cover $B_R$ by parallelepiped $\sq$ of dimensions $K^{-1}R \times \cdots \times K^{-1}R  \times R$ with the long axis perpendicular to $\tau$. We denote the collection of $\sq$ by $\bB_\tau$ and $\bB = \cup_\tau \bB_\tau$. In addition, given $\sq \in  \bB$, we denote by $\tau(\sq)$ the $\tau$ for which $\sq\in \bB_\tau$. 

Let $\T_\sq$ denote the collection of all $T_1 \in \cup_{B\in \Bn} \T(K^2;B)$ such that $2T_1$ intersects $\sq$ and  $\tau(T_1) = \tau(\sq)$. Let $Y_\sq =\cup_{T_1 \in \T_{\sq} } 2T_1$. Let $\W_\sq$ denote the collection of $T\in \W_{\tau(\sq)}$ such that $2T$ intersects $Y_\sq$ and define $f_\sq = \sum_{T\in \W_\sq} f_T$. We record here that \[ \sum_{\sq} \sum_{T\in \W_\sq} \| f_T \|_{L^p}^p \les \sum_{T\in \W} \| f_T \|_{L^p}^p.\] 
Moreover, $f_\sq$ is microlocalized to $(\sq,\tau(\sq))$ and we have
\begin{align*}
 \sum_{B\in \Bn} \sum_{T_1 \in \T(K^2;B)} \|  f_{T_1} \|_{L^p}^p 
&\les \sum_{\sq} \sum_{T_1\in \T_\sq} \| f_{T_1} \|_{L^p}^p \\
& \les \sum_{\sq} \| f_{\tau(\sq)} \|_{L^p(Y_\sq)}^p 
\les \sum_{\sq} \| f_{\sq} \|_{L^p(Y_\sq)}^p. 
\end{align*}

 Summing over $p$-th power of \eqref{eqn:mainb} over narrow $B\in \Bn$, we get
\begin{equation}\label{eqn:decK}
	\begin{split}
		\| f\|_{L^p(Y)}  & \less K^{\eps_1} (\eta K^{-(d+1)})^{\frac{1}{p}-\frac{1}{p_{d-1}}}   K^{(d-2)(\frac{1}{2}-\frac{1}{p})}
		\left(   \sum_{\sq} \|  f_{\sq} \|_{L^p(Y_\sq)}^p \right)^{1/p}.
	\end{split}
\end{equation}  

Next, we apply the induction hypothesis to $f_{\sq}$ after a parabolic rescaling. Without loss of generality, suppose that $\tau(\sq)$ is contained in  $B_{K^{-1}}^{d-1}(0) \times B^1_{K^{-2}}(0)$. Let $L(x_1,\cdots,x_d) = (Kx_1,\cdots,Kx_{d-1},K^2 x_d)$. Then we may apply the induction hypothesis to $f_{\sq}\circ L$ over $L^{-1}(Y_\sq) \subset L^{-1}(\sq)$ at the scale $R_1= R/K^2$:
\begin{equation}\label{eqn:induc}
	\|  f_{\sq} \|_{L^p(Y_\sq)} \les R_1^{\eps/2} (|L^{-1}(Y_\sq)| R_1^{-d})^{\alpha(p)} (R_1)^{\frac{d-1}{2} (\frac{1}{2}-\frac{1}{p})}   \left( \sum_{T\in \W_\sq} \| f_T \|_{L^p}^p \right)^{1/p}. 
\end{equation}

For every $\sq$, 
\[ |\T_\sq| \eta \les \sum_{T_1\in \T_\sq} |3T_1\cap Y| \les |Y|.\]
Thus, we have 
\[ |L^{-1}(Y_\sq)| = K^{-(d+1)} |Y_\sq| \sim K^{-(d+1)}  |T_1| |\T_\sq|  \les K^{O(\eps_1)}  \eta^{-1} |Y|. \] 
By combining this estimate, \eqref{eqn:decK} and \eqref{eqn:induc}, we obtain
\[ 
	\begin{split}
		\| f\|_{L^p(Y)}  \less &R^{\eps/2} K^{O(\eps_1)}  (\eta K^{-(d+1)})^{\frac{1}{p}-\frac{1}{p_{d-1}}}   K^{(d-2)(\frac{1}{2}-\frac{1}{p})}   (\eta^{-1} K^{2d})^{\alpha(p)}  K^{-(d-1) (\frac{1}{2}-\frac{1}{p})}  \\
		&(|Y|R^{-d})^{\alpha(p)} R^{\frac{d-1}{2} (\frac{1}{2}-\frac{1}{p})}  	 \left( \sum_{T\in \W} \| f_T \|_{L^p}^p \right)^{1/p} \\
		\less &R^{\eps/2} K^{O(\eps_1)}  (\eta K^{-(d+1)})^{\frac{1}{p}-\frac{1}{p_{d-1}} - \alpha(p) }   K^{(d-1)\alpha(p) - (\frac{1}{2}-\frac{1}{p})}   \\
		&(|Y|R^{-d})^{\alpha(p)} R^{\frac{d-1}{2} (\frac{1}{2}-\frac{1}{p})}  	 \left( \sum_{T\in \W} \| f_T \|_{L^p}^p \right)^{1/p}.
	\end{split}
\]

Note that 
\[ (\eta K^{-(d+1)})^{\frac{1}{p}-\frac{1}{p_{d-1}} - \alpha(p) } =O(K^{O(\eps_1)}) \]
since $\eta\les |T_1|$ and $\frac{1}{p}-\frac{1}{p_{d-1}} - \alpha(p) \geq 0$.  Moreover, $(d-1)\alpha(p) - (\frac{1}{2}-\frac{1}{p}) \leq 0$ by the choice of $\alpha(p)$. Therefore, we have
\[ 	\| f\|_{L^p(Y)}  \les R^\eps (|Y|R^{-d})^{\alpha(p)} R^{\frac{d-1}{2} (\frac{1}{2}-\frac{1}{p})}  	 \left( \sum_{T\in \W} \| f_T \|_{L^p}^p \right)^{1/p} \]
and the induction closes.

\section{Convergence of Bochner-Riesz means: Proof of \Cref{thm:BR}}\label{sec:BR}
In this section, we prove \Cref{thm:BR}. Given  \Cref{thm:main}, it essentially follows from  the argument of Gan and Wu \cite{GanWu}.  The main goal of this section is to give an exposition of their argument. We fix $j\in \N$.
\begin{defn} A multiplier $m_j$ is type $j$ if 
\begin{equation}\label{eqn:frequency}
 | \partial^\alpha m_j(\xi)| \les_{\alpha,N}2^{j|\alpha|} (1+2^j|1-|\xi||)^{-N}
\end{equation}
and 
\[ \widehat{m_j}(x) = 2^{-j\frac{d+1}{2}} e^{\pm 2\pi i |x| } a(2^{-j}x) \] for some $a \in C_0^\infty$ supported on the annulus $|x|\sim 1$. 
\end{defn}

We use the notation $\widehat{m(D)f}(\xi) =  m(\xi)\widehat{f}(\xi)$.
The Bochner-Riesz means $S_t^\lambda f$ can be decomposed as 
\[  S_t^\lambda f = m_0(D/t) f + \sum_{j=1}^\infty  2^{-\lambda j } m_j(D/t) f, \]
where  $m_j$ is type $j$ for $j\geq 1$ and $m_0(D/t) f$ is dominated by the Hardy-Littlewood maximal function uniformly in $t$ (see  \cite{Stein,TaoW}). Tao \cite{TaoW} reduced the weak-type estimate \eqref{eqn:maxBR}  to a bound on the maximal function $f \mapsto \sup_{t\sim 1} |m_j(D/t)f|$. To be specific, if
\begin{equation}\label{eqn:maxBRreduction1}
	\| \sup_{t\sim 1 } |m_j(D/t) f| \|_{L^p(B(0,C2^j))} \less 2^{j \lambda_0}  \|f\|_{L^p}
\end{equation}
holds some $\lambda_0$ for type $j$ multipliers $m_j$ and  test functions $f$  for all $j\geq 1$, then the weak-type estimate \eqref{eqn:maxBR} holds for $\lambda>\lambda_0$. Here and in the following, we denote by $A\less B$ the inequality $A\leq C_\eps 2^{\eps j} B$ or $A\leq C_\eps R^{\eps} B$ which holds for any $\eps>0$. 

\subsection{Linearization}
In this subsection, we reduce \Cref{thm:BR} to the following.
\begin{prop}\label{prop:BR}
		Fix a collection of disjoint intervals $\{I_1, \cdots, I_{2^j} \}$ of length $\sim 2^{-j}$ for which $ \cup_{l} I_l =[1/2,1]$. Let $c_l \in I_l$ and let $\{ F_1, \cdots, F_{2^j}\}$ be a partition of $B(0,C2^j)$. For a type $j$ multiplier $m_j$, let $T$ denote the linear operator defined as
		\[T f(x) = \sum_{l=1}^{2^j} 1_{F_l}(x) m_j(D/c_l) f(x). \]
		 Then 
	\[ \| T f\|_{L^{\frac{86}{57}}} \less 2^{ \frac{9}{86} j} \| f\|_{L^\frac{86}{57}}, \] 
		where the implicit constant is independent of the choice of $\{c_l, F_l\}$. 
\end{prop}
We postpone the proof of \Cref{prop:BR} to following subsections and first look at its implications. Clearly, $T$  in \Cref{prop:BR} is dominated by the maximal function $f \mapsto \sup_{t\sim 1} |m_j(D/t)f|$. Conversely, we have
\begin{lemma}\label{lem:red}
	Let $T$ be as in \Cref{prop:BR}. Assume that  for any type $j$ multiplier $m_j$, there exists a constant $A_j>0$ independent of $\{F_l, c_l\}_{1\leq l\leq 2^j}$ such that 
	\[  \| T f\|_{L^p} \leq A_j  \| f\|_{L^p}. \]
	Then for every type $j$ multiplier $m_j$,
	\[ \|  \sup_{t\sim 1} |m_j(D/t) f|\|_{L^p(B(0,C2^j))} \les A_j  \| f\|_{L^p}. \]
\end{lemma}
\begin{proof}
Let $m_j$ be a	type $j$ multiplier.	We linearlize and discretize the maximal operator. Fix a measurable function $t:B(0,CR) \to [1/2,1]$. For each $l=1,\cdots,2^j$, define $F_l = \{ x\in B(0,CR)  : t(x) \in I_l \}$ so that 
	\[|m_j(D/t(x))f(x)|^p \leq \sum_{l} 1_{F_l}(x) \sup_{t\in I_l}| m_j(D/t)f(x)|^p.\]
	Let $a_l$ denote the left endpoint of the interval $I_l$. By the fundamental theorem of calculus and H\"{o}lder, we have
	\[  \sup_{t\in I_l} |m_j(D/t)f(x)|^p \les |m_j(D/a_l)f(x)|^p +  |I_l|^{p-1} \int_{I_l} | \partial_t [m_j(D/t)f(x)]|^p dt.\]
	
	For each $1\leq k\leq d$, let $m_{j,k} (\xi) = |I_l| \xi_k \partial_{k} m_{j}(\xi)$. Then $m_{j,k}$ is a type $j$ multiplier and 
	\[  \sup_{t\in I_l} |m_j(D/t)f(x)|^p \les |m_j(D/a_l)f(x)|^p +  |I_l|^{-1} \int_{I_l}\sum_{1\leq k\leq d}  | m_{j,k}(D/t)f(x)|^p dt.\]
	Choose $b_l\in I_l$ for which 
	\[  \int_{F_l} \sum_{1\leq k\leq d} | m_{j,k}(D/b_l)f|^p   \gtrsim  \sup_{t\in I_l} \int_{F_l} \sum_{1\leq k\leq d} | m_{j,k}(D/t)f|^p . \]
	Then \[  \int_{F_l}\sup_{t\in I_l} |m_j(D/t)f|^p   \les \int_{F_l} |m_j(D/a_l)f|^p +  \int_{F_l} \sum_{1\leq k\leq d} | m_{j,k}(D/b_l)f|^p .\]
	By summing the inequality over $l$ and using the assumption,   we get
	\[ \int_{B(0,C2^j)} |m_j(D/t(x)) f|^p  \les A_j^p  \| f\|_{L^p}^p. \]
	 Since this holds for any measurable function $t$, the claim is proved.
\end{proof}

By \Cref{lem:red} and Tao's reduction discussed earlier, we have reduced \Cref{thm:BR} to \Cref{prop:BR}, which will be proved in the following subsections.

\subsection{Decompositions} In this subsection, we work on $\R^d$ for  $d\geq 2$. 
We fix $\{ c_l,F_l\}_{l=1}^{2^j}$  as in \Cref{prop:BR} and let
\[ T f(x) = \sum_{l} 1_{F_l}(x) m_j(D/c_l) f(x). \]
Given any exponent $1\leq p_0\leq \infty$, by dyadic pigeonholing, there exists $\mathcal{J} \subset \{1,2,\cdots, 2^j\}$ and $2^{-100d j} \leq \gamma\leq 2^{jd}$ for which $|F_l|\sim \gamma$ for all $l\in \mathcal{J}$ \footnote{Those $F_l$ with $|F_l| \leq 2^{-100dj}$ is ignored here as these can be handled by crude estimates.}   and
\begin{equation}\label{eqn:pigeonh}
	\| T f\|_{L^{p_0}}\less \| T_{\mathcal{J}} f\|_{L^{p_0}}, 
\end{equation}
where  
\[ T_{\mathcal{J}} f(x) = \sum_{l\in \mathcal{J}} 1_{F_l}(x) m_j(D/c_l) f(x). \]

By real interpolation, it suffices to work with $f=1_E$, the characteristic function of $E\subset B(0,2^j)$. Given $0<\eps\ll 1$, let $0<\eps_1 \ll  \eps$.  For a technical reason, we work with the parameter $R=2^{(1-\eps_1)j}$, but one may identify $R$ and $2^j$ as $R\approx 2^j$. The main result of this subsection is the following. 
\begin{prop}[cf. \cite{GanWu}]\label{prop:mainBR}
	Let $E\subset B(0,2^j)$ and $2\leq p \leq p_d$. Let $R=2^{(1-\eps_1)j}$. Suppose that \eqref{eqn:goal} holds with the exponent $\alpha(p)$. Then 
	\begin{equation}\label{eqn:bound1}
		|\{x:  |T_{\mathcal{J}}1_{E}(x)|  > \alpha\}|  \less \alpha^{-2}  R^{\frac{1}{2p-3} (\frac{d+1}{4}p - \frac{d-1}{2}) -1}  (\gamma R^{-d})^{\frac{p}{2p-3}\alpha(p)} |E|^{\frac{3p-5}{2p-3}}.
	\end{equation}
\end{prop}

\begin{remark}
    When $d=2$ and $p=6$,  \Cref{prop:mainBR} gives the $L^{18/13}$ bound obtained in \cite{LiWu}. However, when $d\geq 3$ and $p=p_d$,   \Cref{prop:mainBR} gives sufficient conditions for the  weak-type estimate \eqref{eqn:maxBR} which are no better than the classical sufficient condition $\lambda >(d-1)(\frac{1}{p}-\frac{1}{2})$. 
\end{remark}
\Cref{prop:mainBR} follows from the argument from \cite{GanWu}. We sketch the proof for the sake of exposition. We decompose $E$ according to the density.  For a parameter $\beta_2>0$ to be chosen, let $\{ Q\}$ denote the collection of maximal dyadic cubes such that 
\[ |E\cap Q| \geq \beta_2 l(Q), \]
where $l(Q)$ is the side-length of $Q$. Let $E_2= \cup_Q (E\cap Q)$ and $E_1 = E \setminus E_2$. Then 
\[ |E_1 \cap B| \leq \beta_2 l(B)\]
for any dyadic cube $B$; if otherwise, there is a dyadic cube $Q$ in the collection containing $B$ leading to the contradiction $|E_1 \cap B|=0$. These sets $E_1$ and $E_2$ are the low and the high density parts of $E$, respectively. Note that $T_{\mathcal{J}}1_E = T_{\mathcal{J}} 1_{E_1} + T_{\mathcal{J}} 1_{E_2}$. 

For the high-density part $E_2$, we have the following estimate from a local $L^2$-estimate, which goes back to \cite{TaoM}.
\begin{lemma}[{\cite[Lemma 3.4]{GanWu}}]\label{lem:L2}
	For any $\alpha>0$, 
	\[ | \{ x :  |T_{\mathcal{J}} 1_{E_2}(x)| \geq \alpha \} \les \alpha^{-2} 2^{-j} \beta_2^{-1} |E|^2.\]	
\end{lemma}

For the low-density part $E_1$, we use a wave-packet decomposition. Let $R=2^{j(1-\eps_1)}$, so that $m_j$ is essentially supported on the $R^{-1}$-neighborhood of $\S^{d-1}$. Cover $\S^{d-1}$ by finitely overlapping caps $\{ \theta\}$ of diameter $R^{-1/2}$ and let $\{ \chi_{\theta} \}$ be an associated smooth partition of unity. Define $m_{j,\theta}(\xi) = m_j(\xi) \chi_{\theta}(\xi/|\xi|)$ so that $m_{j,\theta}(\xi/c_l)$ is essentially a bump function supported on a box $\theta_l$ of dimensions $R^{-1/2}\times \cdots R^{-1/2}\times R^{-1}$. 

Let $f_{\theta_l}$ denote $m_{j,\theta}(D/c_l) f$.
Let $\T_\theta$ denote the collection of tubes as in \Cref{sec:weighted}. Then $f_{\theta_l}$ is morally constant on each $T\in \T_{\theta}$. Indeed, we may choose  an $L^1$-normalized non-negative function $\Psi_\theta$ such that $\Psi_\theta(x) \sim \Psi_\theta(y)$ for all $x,y\in T$ and 
\[ |f_{\theta_l} (x)| \les 2^{O(\eps_1 j)} |f_{\theta_l}|* \Psi_\theta (x) \]
up to a negligible error term (cf. \Cref{sec:broad}).
Therefore $ |f_{\theta_l}|* \Psi_\theta (x) \sim |f_{\theta_l}|* \Psi_\theta (y) $ for all $x,y \in T$.  Let $f=\chi_{E_1}$. Given $\beta_1>0$, we let 
\[ \T_{\theta_l,small} =\{ T \in \T_{\theta} : \| |f_{\theta_l}|* \Psi_\theta \|_{L^\infty(T)}  < \beta_1 \} \]
and $\T_{\theta_l,large} =  \T_{\theta} \setminus \T_{\theta_l,small}$.  Let $\{ \chi_{T} \}_{T\in \T_\theta}$ denote a smooth partition of unity associated with the covering $\T_\theta$.  Note that $f_{\theta_l} \chi_{T}$ is microlocalized to $(T, \theta_l)$ for each $T\in \T_\theta$.  We define
 \[ T_{small}(f) = \sum_{l\in \mathcal{J}} 1_{F_l}  \sum_{\theta} \sum_{T\in \T_{\theta_l,small}} f_{\theta_l} \chi_{T}. \]
We  define $T_{large}(f)$ similarly, which  gives the  decomposition $T_{\mathcal{J}} f = T_{small}(f)+T_{large}(f)$.

\begin{lemma}[{cf. \cite[Lemmas 3.2 and 3.3]{GanWu}}] \label{lem:32} 
	Suppose that  the weighted $\ell^p$-decoupling estimate \eqref{eqn:goal} holds with the exponent $\alpha(p)$. Then
	\begin{align}
		\int |T_{large} (1_{E_1})|  &\less \beta_1^{-1} \beta_2 R^{-1/2}  |E_1|, \nonumber \\
		\int |T_{small} (1_{E_1})|^p dx &\less \left[(\gamma R^{-d})^{\alpha(p)} R^{\frac{d-1}{2}(\frac{1}{2}-\frac{1}{p})}
		\right]^p \beta_1^{p-2}   |E_1|. \label{eqn:small}
	\end{align}
\end{lemma}
The assumption on the weighted $\ell^p$-decoupling inequality \eqref{eqn:goal} is utilized only for the bound \eqref{eqn:small}. For the sake of completeness, we give the proof. 
\begin{proof}[Proof of \eqref{eqn:small}] 
By the weighted $\ell^p$-decoupling estimate \eqref{eqn:goal},
\begin{align*}
		\int |T_{small} (f)|^p dx &= \sum_{l\in \mathcal{J}} \int_{F_l}  | \sum_{\theta} \sum_{T\in \T_{\theta_l,small}} f_{\theta_l} \chi_{T}|^p dx \\
		&\less \left[  ( \gamma R^{-d})^{\alpha(p)} R^{\frac{d-1}{2} (\frac{1}{2}-\frac{1}{p})} \right]^p\sum_{l\in \mathcal{J}}  \sum_{\theta} \sum_{T\in \T_{\theta_l,small}} \|f_{\theta_l} \chi_{T} \|_{L^p}^p.
\end{align*}
 Let $f= 1_{E_1}$. Note that 
\[  \|f_{\theta_l} \chi_{T} \|_{L^\infty} \less  \||f_{\theta_l}| * \Psi_\theta \|_{L^\infty(T)}  \less   \beta_1. \]
Therefore, 
\begin{align*}
	\sum_{l\in \mathcal{J}}  \sum_{\theta} \sum_{T\in \T_{\theta_l,small}} \|f_{\theta_l} \chi_{T} \|_{L^p}^p \less  \beta_1^{p-2} \sum_{l\in \mathcal{J}}   \sum_{\theta} \sum_{T\in \T_{\theta}} \|f_{\theta_l} \chi_{T} \|_{L^2}^2 \less \beta_1^{p-2} \| f\|_{L^2}^2
\end{align*}
by the $L^2$-orthogonality, giving the claimed bound.
\end{proof}

We may now prove \Cref{prop:mainBR}. 
\begin{proof}[Proof of \Cref{prop:mainBR}]
	By \Cref{lem:32},
	\begin{align*}
		& |\{x:  |T_{\mathcal{J}} 1_{E_1}(x)|  > \alpha\}|  \\
		\leq &  |\{x: |T_{large}(1_{E_1})(x)| > \alpha/2\}| +  |\{x: |T_{small}(1_{E_1})(x)| > \alpha/2\}| \\
		\less & \alpha^{-1} \beta_1^{-1} \beta_2 R^{-1/2}  |E| + \alpha^{-p}  (\gamma R^{-d})^{p \alpha(p)} R^{\frac{d-1}{2}(\frac{p}{2}-1)} \beta_1^{p-2}   |E|. 
	\end{align*}
	We take  $\beta_1$ so that the two terms are equal, which gives
	 \begin{equation}\label{eqn:E1}
		|\{x:  |T_{\mathcal{J}} 1_{E_1}(x)|  > \alpha\}| \less  \alpha^{-2} \beta_2^{\frac{p-2}{p-1}} R^{\frac{d-1}{4}\frac{p-2}{p-1}-\frac{1}{2}+\frac{1}{2(p-1)}} (\gamma R^{-d})^{\frac{p}{p-1} \alpha(p)} |E_1|. 
	\end{equation}
	
	By combining  \Cref{lem:L2} and \eqref{eqn:E1},
\begin{align*}
		& |\{x:  |T_{\mathcal{J}} 1_{E}(x)|  > \alpha\}|  \\
		 \leq &		 |\{x:  |T_{\mathcal{J}} 1_{E_1}(x)|  > \alpha/2\}|  + |\{x:  |T_{\mathcal{J}} 1_{E_2}(x)|  > \alpha/2\}| \\
		\less  &   \alpha^{-2} \beta_2^{\frac{p-2}{p-1}} R^{\frac{d-1}{4}\frac{p-2}{p-1}-\frac{1}{2}+\frac{1}{2(p-1)}} (\gamma R^{-d})^{\frac{p}{p-1} \alpha(p)} |E| + \alpha^{-2} R^{-1} \beta_2^{-1} |E|^2.
	\end{align*}
   Finally,	it suffices to choose $\beta_2$ so that the two terms are equal: 
	\[ \beta_2^{-1} = |E|^{-\frac{p-1}{2p-3}} R^{\frac{1}{2p-3} (\frac{d+1}{4}p - \frac{d-1}{2}) } (\gamma R^{-d})^{\frac{p}{2p-3}\alpha(p)}.\]
\end{proof}

\subsection{Proof of \Cref{prop:BR}}
By applying \Cref{prop:mainBR} and  \Cref{thm:main} with $p=\frac{14}{5}$ and $\alpha(p) = \frac{1}{7}$, we obtain  the restricted weak-type estimate 
\[
	|\{x:  |T_{\mathcal{J}} 1_{E}(x)|  > \alpha\}|  \less \alpha^{-2}  R^{-\frac{5}{13}}  (\gamma R^{-2})^{\frac{2}{13}} |E|^{\frac{17}{13}},
\]
or equivalently 
\[ \| T_{\mathcal{J}} f \|_{L^{2,\infty}} \less R^{-\frac{5}{26} } (\gamma R^{-2})^{\frac{1}{13}} \| f\|_{L^{\frac{26}{17},1}}. \]
By H\"{o}lder's inequality (for Lorentz spaces), we get 
\begin{equation}\label{eqn:bound2d1}
	\| T_{\mathcal{J}} f \|_{L^{\frac{26}{17}}} \less2^{2(\frac{17}{26}-\frac{1}{2})j} 2^{-\frac{5}{26} j} (\gamma R^{-2})^{\frac{1}{13}} \| f\|_{L^{\frac{26}{17},1}} = 2^{\frac{3}{26} j} (\gamma R^{-2})^{\frac{1}{13}} \| f\|_{L^{\frac{26}{17},1}}.
\end{equation}

In order to use the gain of $(\gamma R^{-2})^{\frac{1}{13}}$ in \eqref{eqn:bound2d1}, we interpolate it with a $L^{4/3}$-estimate. We recall the $L^{4/3}$-estimate due to Carleson and Sj\"{o}rin \cite{CaSj} and C\'{o}rdoba \cite{CorBR} 
\[ \| m_j(D) f \|_{L^{4/3}} \less \| f\|_{L^{4/3}}. \]
Since the $L^{4/3}$-operator norm of  $m_j(D/t)$ is independent of $t>0$, we have 
\begin{equation}\label{eqn:BRconj}
	\|T_{\mathcal{J}} f\|_{L^{4/3}} = \left( \sum_{l\in \mathcal{J}} \| m_j(D/c_l) f\|_{L^{4/3}(F_l)}^{4/3} \right)^{3/4} \less |\mathcal{J}|^{3/4} \|f \|_{L^{4/3}}.
\end{equation}
We interpolate \eqref{eqn:bound2d1} and  \eqref{eqn:BRconj} so that we may apply $\gamma R^{-2} |\mathcal{J}|  \less 1$. To be specific, take $\theta = \frac{39}{43}$ so that $\theta\cdot \frac{1}{13} = (1-\theta) \frac{3}{4}$. Then $\theta \frac{17}{26} + (1-\theta)  \frac{3}{4}=\frac{57}{86}$. Real interpolation of \eqref{eqn:bound2d1} and  \eqref{eqn:BRconj} gives 
\[ \| T_{\mathcal{J}} f\|_{L^{\frac{86}{57}}} \less 2^{ \frac{9}{86} j}\| f\|_{L^\frac{86}{57}}, \]
which completes the proof of \Cref{prop:BR} by \eqref{eqn:pigeonh}.

\begin{remark}
	Consider the example $f(x) = e^{2\pi i x_2} \psi(x_1,2^{-j/2} x_2)$ from \cite{TaoW}. In this example, if $x_1\sim x_2 \sim 2^j$, then 
	\[ \sup_{t\sim 1} |m_j(D/t)f(x)|\sim  |m_j(D/t(x))f(x)| \sim 2^{-j}, \]
	where $t(x) = |x|/x_2$. Thus, we may choose $\{ c_l,F_l\}$ so that $\gamma \sim 2^{j}$, $|\mathcal{J}|\sim 2^j$, $\| T_\mathcal{J} f\|_{L^p} \gtrsim 2^{2j(\frac{1}{p}-\frac{1}{2})}$ and $\| f\|_{L^q} \sim 2^{j\frac{1}{2q}}$. In particular,
	\begin{align*}
		\|T_{\mathcal{J}} f\|_{L^{4/3}}  &\gtrsim  |\mathcal{J}|^{1/8} \|f \|_{L^{4/3}}.
	\end{align*}
	This suggests that there may be a significant room for an improvement in the estimate \eqref{eqn:BRconj}. Any improvement to \eqref{eqn:BRconj} would advance our knowledge of the almost everywhere convergence of the Bochner-Riesz means. 
\end{remark}

\subsection{An improved bound for $d= 3$}\label{sec:highdim}

\begin{theorem}\label{thm:BR3d}
	Let $d=3$, $1< p< 2$ and \[ \lambda>\max\left(\frac{174}{85p}-\frac{89}{85},\frac{146}{77}\left(\frac{1}{p}-\frac{1}{2}\right) \right).\] 
	Then the weak-type estimate \eqref{eqn:maxBR} holds. Consequently, $\lim_{t\to \infty} S^\lambda_t f(x) = f(x)$ for almost every $x\in \R^3$ for any $f\in L^p(\R^3)$.
\end{theorem}

\Cref{thm:BR3d} gives a small improvement to the result obtained in \cite{GanWu}; when $p=\frac{3}{2}$, \eqref{eqn:maxBR} holds for $\lambda> \frac{27}{85}=0.317\cdots$, improving the sufficient condition $\lambda> \frac{107}{325}=0.329\cdots$ from \cite{GanWu}.

\begin{proof}
Let $E\subset B(0,2^j)$. By taking $p=\frac{14}{5}$ and $\alpha(p)=\frac{1}{14}$ in \Cref{prop:mainBR}, we get 
\[
|\{x:  |T_{\mathcal{J}} 1_{E}(x)|  > \alpha\}|  \less \alpha^{-2}  R^{-\frac{4}{13}} (\gamma R^{-3})^{\frac{1}{13}}  |E|^{\frac{17}{13}},
\]
or equivalently,
\[ 		\| T_{\mathcal{J}} f\|_{L^{2,\infty}} \less 2^{-\frac{2}{13}j} (\gamma R^{-3})^{\frac{1}{26}}  \| f\|_{L^{\frac{26}{17},1}}.\]
By H\"{o}lder's inequality, we get 
\begin{equation}\label{eqn:bound3d1}
		\| T_{\mathcal{J}} f\|_{L^{\frac{26}{17}}} \less  2^{\frac{4}{13}j} (\gamma R^{-3})^{\frac{1}{26}}  \| f\|_{L^{\frac{26}{17},1}}.
\end{equation}

On the other hand, by the sharp $L^p$ estimate for  the Bochner-Riesz means  at the exponent $p=\frac{13}{9}$ due to Wu \cite{Wu_Bochner} (see also \cite{Guo_Bochner}), there holds 
\[  \| m_j(D) f\|_{L^{\frac{13}{9}}} \less  2^{(3(\frac{9}{13}-\frac{1}{2})-\frac{1}{2})j} \|f \|_{L^{\frac{13}{9}}}.\]
Therefore,
\begin{equation}\label{eqn:BR3d}
	\|T_{\mathcal{J}} f\|_{L^{\frac{13}{9}}} = \left( \sum_{l\in \mathcal{J}} \| m_j(D/c_l) f\|_{L^{\frac{13}{9}}(F_l)}^{\frac{13}{9}} \right)^{\frac{9}{13}} \less   2^{\frac{j}{13}} |\mathcal{J}|^{\frac{9}{13}} \|f \|_{L^{\frac{13}{9}}}.
\end{equation}	

We interpolate \eqref{eqn:bound3d1} and  \eqref{eqn:BR3d} so that we may apply $\gamma |\mathcal{J}| R^{-3} \less 1$. This gives 
\[ \| T_{\mathcal{J}} f\|_{L^{\frac{13\cdot 19}{9\cdot 18}}}  \less  2^{\frac{73}{247}j} \| f\|_{L^\frac{13\cdot 19}{9\cdot 18}}. \]
This implies, by \eqref{eqn:pigeonh} and \Cref{lem:red}, that
\begin{equation}\label{eqn:3d}
	\| \sup_{t\sim 1} |m_j(D/t) f|\|_{L^{\frac{247}{162}}(B(0,C2^j) )} \less 2^{\frac{73}{247}j} \| f\|_{L^\frac{247}{162}}. 
\end{equation}
By interpolation with the standard $L^1$ and $L^2$ estimates, we get for $1\leq p \leq 2$, 
\begin{align*}
		\| \sup_{t\sim 1} |m_j(D/t) f|\|_{L^{p}(B(0,C2^j) )} \less \max(2^{j(\frac{174}{85p}-\frac{89}{85})}, 2^{j\frac{146}{77}(\frac{1}{p}-\frac{1}{2})} ) \| f\|_{L^p}. 
\end{align*}
By the reduction of Tao, this implies \Cref{thm:BR3d}.
\end{proof}
	

\newcommand{\etalchar}[1]{$^{#1}$}
\providecommand{\bysame}{\leavevmode\hbox to3em{\hrulefill}\thinspace}
\providecommand{\MR}{\relax\ifhmode\unskip\space\fi MR }
\providecommand{\MRhref}[2]{%
	\href{http://www.ams.org/mathscinet-getitem?mr=#1}{#2}
}
\providecommand{\href}[2]{#2}

\end{document}